\documentclass[11pt]{article}
\usepackage{amsfonts}
\usepackage{mathrsfs}
\usepackage{cite}
\usepackage{amssymb,amsmath} 
 \allowdisplaybreaks

\catcode`!=11
\let\!int\int \def\int{\displaystyle\!int}
\let\!lim\lim \def\lim{\displaystyle\!lim}
\let\!sum\sum \def\sum{\displaystyle\!sum}
\let\!sup\sup \def\sup{\displaystyle\!sup}
\let\!inf\inf \def\inf{\displaystyle\!inf}
\let\!cap\cap \def\cap{\displaystyle\!cap}
\let\!max\max \def\max{\displaystyle\!max}
\let\!min\min \def\min{\displaystyle\!min}
\let\!frac\frac \def\frac{\displaystyle\!frac}
\catcode`!=12

\let\oldsection\section
\renewcommand\section{\setcounter{equation}{0}\oldsection}

\allowdisplaybreaks
\def\pf{\it{Proof.}\rm\quad}

\def\R{\mathbb{R}}

\newtheorem{thm}{Theorem}[section]

\newtheorem{pro}{Proposition}[section]
\newtheorem{lem}{Lemma}[section]

\begin{document}
\title{\Large\bf THE STABILITY OF STRONG VISCOUS  CONTACT DISCONTINUITY TO A FREE BOUNDARY PROBLEM FOR COMPRESSIBLE NAVIER-STOKES EQUATIONS     }
\author{{Tingting Zheng\footnote{Corresponding author, e-mail: asting16@sohu.com}}\\[2mm]
{\small  Computer and Message Science College,
 Fujian Agriculture
and Forestry University,}\\[2mm]
{\small  Fuzhou 350001, P. R. China }}

\date{}

\maketitle

\noindent{\bf Abstract.} This paper is concerned with nonlinear
stability of strong viscous contact discontinuity to a free boundary
problem for the one-dimensional full compressible Navier-Stokes
equations  in half space $[0,\infty)$.  For the case when the local
stability of the contact discontinuities was first studied by
\cite{X},later generalized by \cite{LX}, local stability of weak
viscous contact discontinuity is well-established by
\cite{HMS,HMX,HXY,HZ,HLM2009}, but for the global stability of the
impermeable gas , fewer strong nonlinear wave stability results have
been obtained excluding zero dissipation \cite{MaSX} or $\gamma\to
1$ gas see \cite{HH}. Our main purpose is to deduce the
corresponding nonlinear stability result  by exploiting the
elementary energy method. We will show in this paper that with a
certain class of big perturbation the global stability result of
strong viscous contact discontinuity to Navier-Stokes equations can
be obtained.
\\[2mm]
\\[2mm]
\section{Introduction}
This paper is concerned with a free boundary problem for a
one-dimensional compressible viscous heat-conducting flow in the
half space $\R_+=[0,\infty)$, which is governed by the following
initial-boundary value problem in Eulerian coordinate
$(\tilde{x},t)$:
\begin{equation}
\left\{
\begin{array}{lll}
\tilde\rho_t+(\tilde\rho \tilde u)_{\tilde x}=0,\quad (\tilde x,t)\in\R_+\times\R_+,\\[2mm]
(\tilde\rho \tilde u)_t+(\tilde\rho \tilde u^2+\tilde p)_{\tilde{x}}=\mu \tilde u_{\tilde{x}\tilde{x}}, \\[2mm]
\left(\tilde\rho\left(\tilde e+\frac{\tilde
u^2}{2}\right)\right)_t+\left(\tilde\rho \tilde u\left(\tilde
e+\frac{\tilde u^2}{2}\right)+\tilde p\tilde
u\right)_{\tilde{x}}=\kappa\tilde
\theta_{\tilde{x}\tilde{x}}+(\mu\tilde u\tilde
u_{\tilde{x}})_{\tilde{x}},
\end{array}\right.
\label{1.1}
\end{equation}
where $\tilde\rho$, $\tilde u$ and $\tilde\theta$ are the density,
the velocity and the absolute temperature, respectively, while
$\mu>0$ is the viscosity coefficient and $\kappa>0$ is the
heat-conductivity coefficients, respectively. The pressure $p=\tilde
p(\tilde \rho,\tilde \theta)$ are related by the second law of
thermodynamics. To simplify our problem, we focus our attention on
the perfect gas. In this situation
\begin{eqnarray*}
&&\tilde p(\tilde \rho,\tilde \theta)=R\tilde \theta\tilde\rho,\\
&&\tilde e(\tilde \rho,\tilde \theta)=\frac{R}{\gamma-1}\tilde
\theta +const,
\end{eqnarray*}
where $R>0$ is the gas constant and $\gamma>1$ is the adiabatic
exponent. We consider the system (\ref{1.1}) in the part $\tilde
x>\tilde x(t)$, where $\tilde x=\tilde x(t)$ is a free boundary,
with the following boundary condition
\begin{equation}\label{1.2}
\frac{d\tilde x(t)}{dt}=\tilde u(\tilde x(t),t),\ \tilde x(0)=0,\
\tilde \theta(\tilde x(t),t)=\theta_->0,
\end{equation}
and
\begin{equation}\label{1.3}
(\tilde p-\mu\tilde u_{\tilde x})|_{\tilde x=\tilde x(t)}=p_-,
\end{equation}
which means the gas is attached at the free boundary $\tilde
x=\tilde x(t)$ to the atmosphere with pressure $p_-$ and the initial
data
\begin{equation}\label{1.4}
(\tilde \rho,\tilde u,\tilde \theta)(\tilde x,0)= ({\tilde
\rho}_0,{\tilde u}_0,{\tilde \theta}_0)(\tilde x),\ \lim_{\tilde
x\to +\infty}({\tilde \rho}_0,{\tilde u}_0,{\tilde \theta}_0)(\tilde
x)=(\rho_+,0,\theta_+),
\end{equation}
where $\rho_+,$ $\theta_+$ are  positive constants and
$\theta_0(0)=\theta_-$. Because here we only consider the case of a
single contact discontinuity, we require
\begin{equation}\label{1.5}
p_-=p_+=R\theta_+\rho_+.
\end{equation}

Since the boundary condition (\ref{1.3}) means the particles always
stay on the free boundary $\tilde x=\tilde x(t)$, if we use
$Lagrangian $ coordinates, then the free boundary becomes a fixed
boundary. That is

\begin{equation}
\left\{
\begin{array}{ll}
v_t-u_x=0,\quad (x,t)\in\R_+\times\R_+,\\[2mm]
u_t+\left(\frac{R\theta}{v}\right)_x=\mu\left(\frac{u_x}{v}\right)_x,
\\[4mm]
\frac{R}{\gamma-1}\theta_t+R\frac{\theta}{v}u_x=\kappa\left(\frac{\theta_x}{v}\right)_x+\mu\frac{u_x^2}{v},\\[4mm]
\theta|_{x=0}=\theta_-,\quad t>0,\\[2mm]
\left(\frac{R\theta_-}{v}-\mu\frac{u_x}{v}\right)(0,t)=p_+,\quad t>0,\\[2mm]
(v,u,\theta)|_{t=0}=(v_0,u_0,\theta_0)\to(v_+,0,\theta_+)\quad\mbox{as}\quad
x\to+\infty,
\end{array}
\right.\label{1.6}
\end{equation}
where $v_+$ and $\theta_\pm$ are given positive constants, and
$v_0,\ \theta_0>0$. In fact $v=1/\rho(x,t),\ u=u(x,t),\
\theta=\theta(x,t)$ and $R\theta/v=p(v,\theta)$ are the specific
volume, velocity , temperature and pressure as in (\ref{1.1}).

 In terms of various boundary values,  Matsumura
\cite{Ma2001} classified all possible large-time behaviors of the
solutions for the one-dimensional (isentropic)compressible
Navier-Stokes equations. In the case that $u(0,t)=0$ (resp.
$u(0,t)<0$), the problem is called the {\it impermeable wall} (resp.
{\it outflow}) problem in which the boundary condition of density
can't be imposed.  There have been a lot of works on the asymptotic
behaviors of solutions to the initial-boundary value (or Cauchy)
problem for the Navier-Stokes equations toward these basic waves or
their viscous versions, see, for example, [3--25] and the references
therein.

   On the other hand, the problem of stability of
contact discontinuities are associated with linear degenerate fields
and are less stable than the nonlinear waves for the inviscid system
(Euler equations). It was observed in \cite{X,LX}, where the
metastability of contact waves was studied for viscous conservation
laws with artificial viscosity, that the contact discontinuity
cannot be the asymptotic state for the viscous system, and a
diffusive wave, which approximated the contact discontinuity on any
finite time interval, actually dominates the large-time behavior of
solutions. The nonlinear stability of contact discontinuity for the
(full) compressible Navier-Stokes equations was then investigated in
\cite{HMS,HZ}
 for the free  boundary value problem and \cite{HMX,HXY}for the Cauchy problem.

As it is shown in the references, we construct viscous contact wave
with artificial viscousity by the
 corresponding Euler system of (\ref{1.6}) with Riemann initial
data which reads as follows:
\begin{equation}
\left\{
\begin{array}{ll}
v_t-u_x=0,\\[2mm]
u_t+p(v,\theta)_x=0,\\[2mm]
\frac{R}{\gamma-1}\theta_t+R\frac{\theta}{v}u_x=0,\\[4mm]
(v,u,\theta)(x,0)=(v_-,0,\theta_-)\quad\mbox{if}\quad x<0,\\[2mm]
(v,u,\theta)(x,0)=(v_+,0,\theta_+)\quad\mbox{if}\quad x>0.
\end{array}
\right.\label{1.7}
\end{equation}

 Because the corresponding Euler
equations (\ref{1.7}) with the Riemann initial data has the
following soluitons
\begin{equation}
\left(\overline{V},\overline{U},\overline{\Theta}\right) =\left\{
\begin{array}{ll}
(v_-,0,\theta_-),\quad x<0,\\[2mm]
(v_+,0,\theta_+),\quad x>0,
\end{array}
\right. \label{1.8}
\end{equation}
provided that
\begin{equation}
p_-=R\frac{\theta_-}{v_-}=p_+=R\frac{\theta_+}{v_+},\label{1.9}
\end{equation}
as that in \cite{HMS} we conjecture that the asymptotic limit
$(V,U,\Theta)$ of (\ref{1.6}) is as follows

\begin{equation}
P(V,\Theta)=R\frac{\Theta}{V}=p_+,\;\;\;U(x,t)=\frac{\kappa(\gamma-1)\Theta_x}{\gamma
R\Theta}, \label{1.10}
\end{equation}
and $\Theta $ is the solution of the following  problem
\begin{equation}
\left\{ \begin{array}{ll}
\Theta_t=a(\ln\Theta)_{xx},\quad a=\frac{\kappa p_+(\gamma-1)}{\gamma R^2}>0,\\[2mm]
\Theta(0,t)=\theta_-,\\[2mm]
 \Theta(x,0)=\Theta_{0}\to \theta_+,\ \ \mathrm{as}\ x\to +\infty
\end{array}\right.\label{1.11}
\end{equation}
with $\Theta_0=\theta_+-(\theta_+-\theta_-)\exp\{1-(1+\alpha
x)^{\delta_0}\}$. It is easy to check that there exist positive
constant  $M_0$ which is independent of $\delta_0$ and $\alpha$ such
that
\begin{eqnarray}\label{1.12b}
&&\|\Theta_0-\theta_+\|_{L^1}\leq M_0\alpha^{-1}\sum_{n=0}^{[\frac{1}{\delta_0}]-1}\prod_{i=0}^n(\frac{1}{\delta_0}-i),\nonumber\\
&&|\Theta_{0x}|\leq M_0\alpha\delta_0,\nonumber\\
&&|\Theta_{0xx}|\leq M_0\alpha^2\delta_0,\nonumber\\
&&\|\Theta_{0x}\|^2\leq M_0\alpha\delta_0,\nonumber\\
&&\|\Theta_{0x}\|_{L^1(\R_+)}\leq M_0,\nonumber\\
&&\|\Theta_{0xx}\|^2+\|(\ln\Theta_0)_{xx}\|^2\leq
M_0\alpha^3\delta_0^2,\nonumber\\
&&\|\Theta_{0xxx}\|^2+\|(\ln\Theta_0)_{xxx}\|^2\leq
M_0\alpha ^5,\nonumber\\
&&\int_{\R_+}\Theta_{0x}^2(1+\alpha x)dx\leq M_0\alpha\delta_0,
\end{eqnarray}
here and following $\alpha$  is a positive constant which will be
determined  in Lemma \ref{lem2.4} due to the artificial viscosity
and $\delta_0$ is a small positive constant which is independent of
$\theta_{\pm}$. To sum up, we have constructed a pair of functions
$(V,U,\Theta)$ such that
\begin{equation}
\left\{
\begin{array}{ll}
R\frac{\Theta}{V}=p_+,\\[4mm]
V_t=U_x,\\[2mm]
U_t+P(V,\Theta)_x=\mu\left(\frac{U_x}{V}\right)_x+F,\\[4mm]
\frac{R}{\gamma-1}\Theta_t+R\frac{\Theta}{V}U_x=\kappa\left(\frac{\Theta_x}{V}\right)_x+\mu\frac{U_x^2}{V}+G,\\[4mm]
(V,U,\Theta)(0,t)=(v_-,\frac{\kappa(\gamma-1)}{\gamma
R}\frac{\Theta_{x}}{\Theta}|_{x=0},\theta_-),\\[4mm]
(V,U,\Theta)(x,0)=(V_0,U_0,\Theta_0)=(\frac{R}{p_+}\Theta_0,\frac{\kappa(\gamma-1)}{\gamma
R}\frac{\Theta_{0x}}{\Theta_0},\Theta_0)\to (v_+,0,\theta_+),\ as\ \
x\to+\infty,
\end{array}\right.\label{1.12}
\end{equation}
where
\begin{eqnarray}G&=&-\mu\frac{U_x^2}{V}=O((\ln\Theta)^2_{xx}),\nonumber\\
F(x,t)&=&\frac{\kappa(\gamma-1)}{\gamma
R}\left\{(\ln\Theta)_{xt}-\mu\left(\frac{(\ln\Theta)_{xx}}{V}\right)_x\right\}\nonumber\\
&=&\frac{\kappa a(\gamma-1)-\mu p_+\gamma}{R
\gamma}\left(\frac{(\ln\Theta)_{xx}}{\Theta}\right)_x.\label{1.13}\end{eqnarray}
We shall show in the next section that $(V,U,\Theta)$ approximates
$(\overline{V},\overline{U},\overline{\Theta})$ in $L^p$ norm with
$p\geq 1$ on any finite time interval as the heat conductivity
$\kappa$ goes to zero. So, we call $(V,U,\Theta)$ the viscous
contact wave for the Navier-Stokes system (\ref{1.6}). The
definition can be more precise according to whether $H(\R_+)$--norm
of the initial perturbation $(\varphi_0(x),\psi_0(x),\zeta_0(x))$
and (or) $|\theta_+-\theta_-|$ big or not, the stability results are
classified into global (or local) stability of strong (or weak )
viscous contact wave.

Our main purpose is to justify that the solution $(v,u,\theta)$ of
the Navier-Stokes system (\ref{1.3}) asymptotically tends to the
strong viscous contact discontinuity $(V,U,\Theta)$. Roughly
speaking, the main result is :``{\it if the oscillation of
temperature and density are not small, the viscous contact
discontinuity is asymptotic stable}''. To deduce the desired
nonlinear stability result by the elementary energy method as in
\cite{HMS,HMX,HXY,HZ,HLM2009,HH,TJJ2011,TZ2012}, it is sufficient to
deduce certain uniform (with respect to the time variable $t$)
energy type estimates on the solution $(v(x,t),u(x,t),\theta(x,t))$
and how to  establish the Poincar$\acute{e}$ type inequality in
Lemma \ref{lem3.2}  without the smallness of $|\theta_+-\theta_-|$
which the arguments employed in
\cite{HMS,HMX,HXY,HZ,HLM2009,TJJ2011,TZ2012} is to use both
smallness $|\theta_+-\theta_-|$  and $N(t)=\sup_{0\leq \tau\leq
t}\|(\varphi,\psi,\zeta)\|_{H^1}$ to overcome such  difficulties.
One of the key points in such an argument is that, based on the  a
priori assumption that $\sup_{0\leq\tau\leq
t}\|(\varphi,\psi,\zeta)\|_{H^1}(\tau)$ is sufficiently small, one
can deduce a uniform lower and upper positive bounds on the specific
volume $v(x,t)$ and temperature $\theta(x,t)$. With such a bound on
$v$ and $\theta$ in hand, one can deduce  a priori $H(\R_+)$ energy
type estimates on $(\varphi,\psi,\zeta)$ in terms of the initial
perturbation $(\varphi_0,\psi_0,\zeta_0)$ provided that
$|\theta_+-\theta_-|$ suitably small, so the stability of weak
contact discontinuity can be obtained . In fact if $N(t)$ not small
and the perturbation of
$\|(\varphi_{0x},\psi_{0x},\zeta_{0x})\|_{L^2(\R_+)}$ not small (see
\cite{HH}), the combination of  the  analysis similar as above with
the standard continuation argument, it also can obtain the upper and
lower bounds of $(v,\theta)$, then that yields the global stability
of strong viscous contact discontinuity for the one-dimensional
compressible Navier-Stokes equations in the condition of $\gamma \to
1$. In all, after researching  the references carefully we find it
is important to get the uniform time estimates of the viscous
contact discontinuity we constructed, then we can obtain the energy
estimates we expected, then the upper and lower bounds of
$(v,\theta)$ come out  .
 So the global stability result can be obtained. It is easy to see that in such a result,
  for all $t\in\R_+$, Osc
 $\theta(t):=\sup_{x\in\R_+}\theta(x,t)-\inf_{x\in\R_+}\theta(x,t)\geq|\theta_+-\theta_-|$,
 the oscillation of the  temperature $\theta(x,t)$ should not be
 sufficiently small when $|\theta_+-\theta_-|$ is not small. Similarly, we can also obtain the oscillation of the  density $\rho(x,t)$ should not be
 sufficiently small too .

Following the above analysis, the rest of this paper is out lined as
follows. In section 2 we study the properties of the  viscous
continuity $(V,U,\Theta)$ in (\ref{1.6}). In section 3, we
reformulate the problem and give the precise statement of our main
theorem. Finally, we complete the proof of the main result by the
global a priori estimates established in section 1.

Throughout this paper, we shall denote $H^l(\R_+)$ the usual $l-th$
order Sobolev space with the norm
$$\|f\|_{l}=\big(\sum_{j=0}^l\|\partial_x^j f\|^2\big)^{1/2},\ \ \|\cdot\|:=\|\cdot\|_{L^2(\R_+)}.$$
For simplicity, we also use $C$ or $C_i$ ($i=1,2,3.....$) to denote
the various positive generic constants.  $C(z)$ stands for constant
about $z$ and $\lim_{z\to 0}C(z)=0$. $\epsilon$ and
$\epsilon_i$($i=1,2,3.....$) stand for a small positive constant in
Cauchy-Schwarz inequality and
$\partial_x^i=\frac{\partial^i}{\partial x^i}$.

\section{Preliminaries}
This section is devoted to study of the viscous contact
discontinuity $(V,U,\Theta)$ in (\ref{1.12}). To finish it, we
construct a parabolic equation about
 $\theta_2$ which play an important role in the time estimates of  $\partial_x^i\Theta$ $(i=1,2,3)$ , it is shown as follows.

\begin{lem}\label{lem2.2}If $\delta_0$ and $\Theta_0$ satisfying the
condition in Theorem \ref{thm1.1} and
\begin{eqnarray*}&&\theta_2(x,t)=\int_{0}^{+\infty}(4\pi
at)^{-1/2}(\Theta_0(h)-\theta_-)\left\{\exp\{-\frac{(h-x)^2}{4at}\}-\exp\{-\frac{(h+x)^2}{4at}\}\right\}\
dh+\theta_-,\end{eqnarray*}
 we can get
 \begin{eqnarray}
 &&\theta_{2t}=a\theta_{2xx};\nonumber\\
 &&\theta_2(0,t)=\theta_-;\nonumber\\
 &&\theta_2(x,0)=\theta_{20}(x)=\left\{
                      \begin{array}{ll}
                        \Theta_{0}(x)\to \theta_+ , & \hbox{$x>0$;} \\
                       -\Theta_{0}(-x)+2\theta_-\to 2\theta_--\theta_+ , & \hbox{$ x\leq 0$,}
                      \end{array}
                    \right.\label{2.2}
 \end{eqnarray}
 and\begin{eqnarray}
 &&\int_0^t\|\theta_{2x}\|^2dt\leq C(1+t)^{1/2},\label{2.3}\end{eqnarray}

\end{lem}
\pf Because $\theta_2(x,t)$ can be rewrite to
$$\theta_2(x,t)=\int_{-\infty}^{+\infty}(4\pi at)^{-1/2}\theta_{20}(h)\exp\{-\frac{(x-h)^2}{4at}\}dh,$$
and $\theta_{20}(x)\in  C^1(\R)$, we find that $\theta_2(x,t)$ is a
fundamental solution of (\ref{2.2}), it is easy to check $\lim_{t\to
0}\theta_2(x,t)=\theta_{20}(x)$, so we finish (\ref{2.2}).

Because \begin{eqnarray}\label{2.5} \theta_{2x}
&=&\int_0^{+\infty}(4\pi
at)^{-1/2}\left(\Theta_0(z)-\theta_-\right)\exp\{-\frac{(z-x)^2}{4at}\}\frac{z-x}{2at}dz\nonumber\\
&&\quad+\int_0^{+\infty}(4\pi
at)^{-1/2}\left(\Theta_0(z)-\theta_-\right)\exp\{-\frac{(z+x)^2}{4at}\}\frac{z+x}{2at}dz\nonumber\\
 &=&\int_0^{\infty}(4\pi
at)^{-1/2}\Theta_{0z}(z)\exp\{\frac{-(z-x)^2}{4at}\}dz\nonumber\\
&&\quad-\int_0^{\infty}(4\pi
at)^{-1/2}\Theta_{0z}(z)\exp\{\frac{-(z+x)^2}{4at}\}dz
,\end{eqnarray}
  we use H$\ddot{o}$lder inequality and Fubini Theorem and
$|\theta_+-\theta_-|\leq\|\Theta_{0z}\|_{L^1(\R_+)}<C$ to
(\ref{2.5}), then we can get
\begin{eqnarray*}\int_0^t\int_0^{\infty}\theta_{2x}^2dxdt&\leq&
C\int_0^t\int_0^{\infty}(4\pi
at)^{-1}\left(\int_0^{\infty}\Theta_{0z}\left(\exp\{-\frac{(z-x)^2}{4at}\}-\exp\{-\frac{(z+x)^2}{4at}\}\right) dz\right)^2dxdt\\
&\leq&C\int_0^t\int_0^{\infty}(4\pi
at)^{-1}\int_0^{\infty}|\Theta_{0z}|\exp\{-\frac{(z-x)^2}{4at}\}
dzdx\int_0^{\infty}|\Theta_{0z}|dzdt\\
&&\quad+C\int_0^t\int_0^{\infty}(4\pi
at)^{-1}\int_0^{\infty}|\Theta_{0z}|\exp\{-\frac{(z+x)^2}{4at}\}
dzdx\int_0^{\infty}|\Theta_{0z}|dzdt\\
 &\leq &C\sqrt{1+t}.\end{eqnarray*}

So we finish this lemma.$\Box$

 Now let's consider the time estimates of $\partial_x^i\Theta$ $(i=1,2,3)$ of
 (\ref{1.11}), we have the following results.
\begin{lem}\label{lem2.3}
If $\Theta_{0x}$ satisfying the condition of (\ref{1.12b}) and a
positive constant $M_0$  is independent of $\delta_0$ and $\alpha$ ,
there exist a positive constant  $C$ such that
\begin{eqnarray}&&\|(\ln\Theta)_x\|^2+a\int_0^t\ \|(\ln\Theta)_{xx}\|^2\ dt \leq
M_0\alpha\delta_0.\label{2.7}\\
&&_{(see (\ref{2.17})-(\ref{2.19}))}\nonumber\\
 &&\|\Theta-\theta_2\|^2+\int_0^t\ \|(\ln\Theta)_{x}\|^2\ dt\leq
C(1+t)^{1/2}.\label{2.8}\\
&&_{{(see (\ref{2.14})-(\ref{2.15}))}}\nonumber\\
 && \|(\ln\Theta)_x\|^2\leq
C(1+t)^{-1/2}.\label{2.9}\\&&_{{(see
(\ref{2.20})-(\ref{2.23}))}}\nonumber\\
&&\|(\ln\Theta)_{xx}\|^2\leq C(1+t)^{-3/2}.\label{2.10}\\&&_{(see
(\ref{2.24})-(\ref{2.28}))}\nonumber\\
&&\|(\ln\Theta)_{xx}\|^2(1+t)+\int_0^t\|\partial^3_x\ln\Theta\|^2(1+t)\
dt\leq
C\delta_0^2.\label{2.11}\\&&_{(see(\ref{2.29}))}\nonumber\\
&& \|\partial^3_x\ln\Theta\|^2\leq
C(1+t)^{-5/2}.\label{2.12}\\&&_{(see
(\ref{2.30})-(\ref{2.32}))}\nonumber\\
&&\int_{\R_+}\Theta_x^2xdx\leq C\delta_0.\label{2.13}\\
&&_{(see (\ref{2.33})-(\ref{2.35}))} \nonumber
\end{eqnarray}
\end{lem}
\pf From (\ref{1.11}) we know
\begin{equation}(\ln\Theta)_t=a\frac{(\ln\Theta)_{xx}}{\Theta},\label{2.16}\end{equation}
both side of it multiply by $(\ln\Theta)_{xx}$ and integrate in
$\R_+\times (0,t)$ we can get

\begin{eqnarray}&&\|(\ln\Theta)_x\|^2+a\int_0^t\ \|(\ln\Theta)_{xx}\|^2\
dt\nonumber\\
&&\leq
C\|(\ln\Theta_{0})_x\|^2+\int_0^t(\ln\Theta)_t(\ln\Theta)_x\big|_0^{\infty}\
dt.\label{2.17}\end{eqnarray}

So we can get
\begin{equation}\|(\ln\Theta)_x\|^2+a\int_0^t\ \|(\ln\Theta)_{xx}\|^2\ dt\leq
M_0\alpha\delta_0.\label{2.19}\end{equation}

Then  if $\int_0^t\int_{\R_+}
\left((\ref{1.11})_1-(\ref{2.2})_1\right)\times(\Theta-\theta_2)dxdt$
  combine with Cauchy-Schwarz inequality we can get
\begin{eqnarray}&&\|\Theta-\theta_2\|^2+\int_0^t\|(\ln\Theta)_x\|^2\
dt\leq C\int_0^t\|{\theta_2}_x\|^2dt.\label{2.14}\end{eqnarray} Use
(\ref{2.3}) to (\ref{2.14}) we can get
\begin{equation}\|\Theta-\theta_2\|^2+\int_0^t\ \|(\ln\Theta)_{x}\|^2\ dt\leq
C(1+t)^{1/2}. \label{2.15}\end{equation} That is (\ref{2.8}).

 Next, from $$\int_0^t\int_{\R_+}(\ref{1.11})_1 \times
\Theta^{-1}(\ln\Theta)_{xx}(1+t)dxdt,$$  we can get
\begin{eqnarray}&&\int_0^t(1+t)\left((\ln\Theta)_t(\ln\Theta)_x\right)(0,t)\ dt\nonumber\\
&&=a\int_0^t\int_0^{\infty}\ \frac{(\ln\Theta)^2_{xx}}{\Theta}(1+t)\
dxdt+\int_0^t\int_0^{\infty}\ \left((\ln\Theta)^2_x\right)_t(1+t)\
dxdt.\label{2.20}\end{eqnarray} Because
\begin{equation}\label{2.21}\int_0^t(1+t)(\ln\Theta)_t(\ln\Theta)_x(0,t)\
dt=0,\end{equation} we can get
\begin{eqnarray}
&&(1+t)\|(\ln\Theta)_x\|^2+\int_0^t\
\int_0^{\infty}(1+t)(\ln\Theta)^2_{xx}\ dx\ dt\nonumber\\
 &&\leq C\|\Theta_{0x}\|^2+\int_0^t\ \int_0^{\infty}\ (\ln\Theta)^2_x\ dx\ dt.\label{2.22}\end{eqnarray}
Combine with (\ref{2.15}) we can get
\begin{eqnarray}
&&(1+t)\|(\ln\Theta)_x\|^2+\int_0^t\
\int_0^{\infty}(1+t)(\ln\Theta)^2_{xx}\ dx\ dt\nonumber\\
 &&\leq C(1+t)^{1/2}.\label{2.23}\end{eqnarray}
That means $\|(\ln\Theta)_x\|^2\leq C(1+t)^{-1/2}$, which
is(\ref{2.9}).

Again from (\ref{1.11})$_1$ we can get
\begin{equation}(\ln\Theta)_{xt}=a\left(\frac{(\ln\Theta)_{xx}}{\Theta}\right)_x.\label{2.24}
\end{equation}
Both side of (\ref{2.24})multiply $\partial^3_x\ln\Theta$ and get
\begin{equation}\left((\ln\Theta)_{xt}\partial^2_x(\ln\Theta)\right)_x-1/2(\partial^2_x\ln\Theta)_t
=a\left(\frac{(\ln\Theta)_{xx}}{\Theta}\right)_x\partial^3_x(\ln\Theta).\label{2.25}\end{equation}
Because
\begin{eqnarray*}&&\left((\ln\Theta)_{xt}\partial^2_x(\ln\Theta)\right)_x(1+t)^2\\
&&=\left((\ln\Theta)_{xt}(\ln\Theta)_{xx}\right)_x(1+t)^2\\
&&=a^{-1}\left((\ln\Theta)_{xt}\Theta_t\right)_x(1+t)^2\end{eqnarray*}
then both side of (\ref{2.25}) multiply $(1+t)^2$ then integrate in
$\R_+\times(0,t)$ and combine with $\Theta_t(0,t)=0$,
$\Theta_t(\infty,t)=0$, $\Theta_x(\infty,t)=0$ and Cauchy-Schwarz
inequality to get for some small $\epsilon>0$ we have
\begin{eqnarray}
&&0\geq a\int_0^t\ \int_0^{\infty}\
\frac{(\ln\Theta)_{xxx}^2}{\Theta}(1+t)^2\ dx\
dt\nonumber\\
&&\quad-\epsilon\int_0^t\ \int_0^{\infty}\
(1+t)^2(\ln\Theta)_{xxx}^2\ dx\ dt-C\epsilon^{-1}a\int_0^t\
\int_0^{\infty}\ (1+t)^2(\ln\Theta)_{xx}^2(\ln\Theta)_x^2\ dx\
dt\nonumber\\
&&\quad+1/2\|(\ln\Theta)_{xx}\|^2(1+t)^2-1/2\|(\ln\Theta_0)_{xx}\|^2-\int_0^t\
\|(\ln\Theta)_{xx}\|^2(1+t)\ dx\nonumber\\
&&\geq Ca\int_0^t\ \int_0^{\infty}\
\frac{(\ln\Theta)_{xxx}^2}{\Theta}(1+t)^2\ dx\
dt\nonumber\\
&&\quad-C\epsilon^{-1}a\int_0^t\ \int_0^{\infty}\
(1+t)^2\|(\ln\Theta)_{xx}\|\|(\ln\Theta)_{xxx}\|(\ln\Theta)_x^2\ dx\
dt\nonumber\\
&&\quad+1/2\|(\ln\Theta)_{xx}\|^2(1+t)^2-1/2\|(\ln\Theta_0)_{xx}\|^2-\int_0^t\
\|(\ln\Theta)_{xx}\|^2(1+t)\ dx.\label{2.26}
\end{eqnarray}
 Take (\ref{2.23})
 into (\ref{2.26}) we can get

\begin{eqnarray}&&\|(\ln\Theta)_{xx}\|^2(1+t)^2+\int_0^t\ \int_0^{\infty}\
(1+t)^2(\ln\Theta)_{xxx}^2\ dx\ dt\nonumber\\
&& \leq C(1+t)^{1/2},\label{2.27}\end{eqnarray} which also means

\begin{equation}\|(\ln\Theta)_{xx}\|^2\leq
C(1+t)^{-3/2},\label{2.28}\end{equation} and finish (\ref{2.10}).

If both side of (\ref{2.25}) multiply by $(1+t)$, similar as the
proof of (\ref{2.27}),  when combine with (\ref{2.19}) we can get
\begin{equation}\|(\ln\Theta)_{xx}\|^2(1+t)+\int_0^t\ \int_0^{\infty}\ (1+t)(\partial^3_x\ln\Theta)^2\ dx\
dt\leq C\delta_0^2,\label{2.29}\end{equation} which means
(\ref{2.11}).

From (\ref{2.24}) we can get
\begin{equation}
\partial_t(\ln\Theta)_{xx}=a\partial^2_x\left(\frac{(\ln\Theta)_{xx}}{\Theta}\right).\label{2.30}
\end{equation}
Because
\begin{eqnarray*}
\left((\ln\Theta)_{xxt}(\ln\Theta)_{xxx}\right)_x=\left(a^{-1}\Theta_{tt}(\ln\Theta)_{xxx}\right)_x
\end{eqnarray*}
when both side of (\ref{2.30}) multiply
$\partial^4_x\ln\Theta(1+\tau)^3$ then integrate in
$\R_+\times(0,t)$ we can get

\begin{eqnarray}
&&\int_0^t\int_0^{\infty}\left(a^{-1}\Theta_{tt}(\ln\Theta)_{xxx}\right)_x(1+\tau)^3dxd\tau\nonumber\\
&&=\int_0^t\int_0^{\infty}a\partial^2_x\left(\frac{(\ln\Theta)_{xx}}{\Theta}\right)\partial^4_x\ln\Theta(1+\tau)^3dxd\tau\nonumber\\
&&\quad+\int_0^t\int_0^{\infty}\frac{1}{2}\left((\partial^3_x\ln\Theta)^2\right)_t(1+\tau)^3dxd\tau.\label{2.31}
\end{eqnarray}

So use (\ref{2.23}) and (\ref{2.27}) we can get that for a small
$\epsilon>0$, (\ref{2.31}) can be change to
\begin{eqnarray*}
&&\|\partial^3_x\ln\Theta\|^2(1+t)^3+C\int_0^t\
(1+\tau)^3\|\partial^4_x\ln\Theta\|^2\ d\tau\\
&&\leq C+C\int_0^t\ \int_0^{\infty}\
(\partial^3_x\ln\Theta)^2(\ln\Theta)_x^2(1+\tau)^3\ dx\
d\tau+C\int_0^t\
\int_0^{\infty}\ (\ln\Theta)_{xx}^4(1+\tau)^3\ dx\ d\tau\\
&&\quad+C\int_0^t\ \int_0^{\infty}\
(\partial^2_x\ln\Theta)^2(\ln\Theta)_x^4(1+\tau)^3\ dx\
d\tau+C\int_0^t\
\int_0^{\infty}\ (\partial_x^3\ln\Theta)^2(1+\tau)^2\ dx\ d\tau\\
&&\leq C\int_0^t\
\|(\ln\Theta)_x\|^2\|\partial^3_x\ln\Theta\|\|\partial^4_x\ln\Theta\|(1+\tau)^3\
d\tau+C\int_0^t\
\|(\ln\Theta)_{xx}\|^3\|\partial^3_x\ln\Theta\|(1+\tau)^3\ d\tau\\
&&\quad+\int_0^t\
\|(\ln\Theta)_{xx}\|^4\|(\ln\Theta)_x\|^2(1+\tau)^3\
d\tau+C(1+t)^{1/2}\\
&&\leq \epsilon\int_0^t\ \|\partial^4_x\ln\Theta\|^2(1+\tau)^3\
d\tau+C\epsilon^{-1}\int_0^t\ \|\partial^3_x\ln\Theta\|^2(1+\tau)^2\ d\tau\\
&&\quad+C\epsilon^{-1} \int_0^t\
\|\partial^2_x\ln\Theta\|^2(1+\tau)\ d\tau+C(1+t)^{1/2}.
\end{eqnarray*}
Again using (\ref{2.23}) and (\ref{2.27}) we can get
\begin{equation}\|\partial^3_x\ln\Theta\|^2(1+t)^3+\int_0^t\
(1+\tau)^3\|\partial^4_x\ln\Theta\|^2\ d\tau\leq
C(1+t)^{1/2}.\label{2.32}\end{equation} This means (\ref{2.12})
finished.

 Now both side of (\ref{2.16}) multiply by
$(\ln\Theta)_{xx}(x-\beta \tau)$ ($\beta>0$) and integrate in
$[\beta\tau,\infty)\times (0,t)$ we can get
\begin{eqnarray}\label{2.33}
&&\int_0^t\int_{\beta\tau}^{\infty}\left((\ln\Theta)_\tau(\ln\Theta)_x(x-\beta\tau)\right)_xdxd\tau-1/2\int_0^t\int_{\beta\tau}^{\infty}\left((\ln\Theta)_x^2(x-\beta\tau)\right)_\tau dxd\tau\nonumber\\
&&\quad-a\int_0^t\int_{\beta\tau}^{\infty}(\ln\Theta)_{xx}(\ln\Theta)_x\Theta^{-1}dxd\tau
-\frac{\beta}{2}\int_0^t\int_{\beta\tau}^{\infty}(\ln\Theta)_x^2dxd\tau\nonumber\\
&&\quad-\int_0^t\int_{\beta\tau}^{\infty}a(\ln\Theta)_{xx}^2(x-\beta\tau)\Theta^{-1}
dxd\tau=\sum_{i=1}^5K_i=0.
\end{eqnarray}
Use Cauchy-Schwarz inequality
\begin{eqnarray*}
|K_3|\leq
\frac{\beta}{4}\int_0^t\int_{\beta\tau}^{\infty}(\ln\Theta)_x^2dxd\tau+C\int_0^t\int_{\beta\tau}^{\infty}(\ln\Theta)_{xx}^2dxd\tau,
\end{eqnarray*}
then combine with (\ref{2.7}) we can get
$$|K_3|+K_4\leq C\int_0^t\int_{\beta\tau}^{\infty}(\ln\Theta)_{xx}^2dxd\tau-\frac{\beta}{4}\int_0^t\int_{\beta\tau}^{\infty}(\ln\Theta)_x^2dxd\tau\leq C\alpha\delta_0-\frac{\beta}{4}\int_0^t\int_{\beta\tau}^{\infty}(\ln\Theta)_x^2dxd\tau.$$
From $K_1$ to $K_5$, (\ref{2.33}) can be change to
\begin{eqnarray}\label{2.34}
&&\int_{\beta t}^{\infty}(\ln\Theta)_x^2(x-\beta t)dx+\int_0^t\int_{\beta\tau}^{\infty}(\ln\Theta)_{xx}^2(x-\beta\tau)dxd\tau\nonumber\\
&&\quad+\frac{\beta}{4}\int_0^t\int_{\beta\tau}^{\infty}(\ln\Theta)_x^2dxd\tau
\leq C\delta_0.
\end{eqnarray}

Because
\begin{eqnarray*}&&\lim_{\beta\to 0}\left|\int_{\beta t}^{\infty}(\ln\Theta)_x^2(x-\beta t)dx-\int_{0}^{\infty}(\ln\Theta)_x^2xdx\right|\\
&=&\lim_{\beta\to 0}\left|\int_{\beta
t}^{\infty}(\ln\Theta)_x^2(x-\beta t-x)dx+\int_{0}^{\beta
t}(\ln\Theta)_x^2xdx\right|=0,
\end{eqnarray*}
we can get \begin{equation}\label{2.35}0<\lim_{\beta\to
0}\int_{\beta t}^{\infty}(\ln\Theta)_x^2(x-\beta
t)dx=\int_{0}^{\infty}(\ln\Theta)_x^2xdx\leq
C\delta_0.\end{equation} So we finish this lemma.$\Box$

The next lemma is concerned with the relation ship between the
viscous continuity and the contact discontinuity. We shall show that
as the heat conductivity $\kappa$ goes to zero, $(V,U,\Theta)$ will
approximate $(\overline{V},\overline{U},\overline{\Theta})$ in
$L^p(\R_+)$ $(p\geq 1)$ norm on any finite time interval.

\begin{lem}\label{lem2.4}
For any given $T\in(0,+\infty)$ independent of $\kappa$ such that
for any $p\geq 1$ and $t\in[0,T]$,
$$\|(V-\overline{V},U-\overline{U},\Theta-\overline{\Theta})\|_{L^p(\R_+)}\to 0,\ \ \mathrm{as}\ \ \kappa\to 0.$$
\end{lem}
\pf  By the definition of $\overline{\Theta}$ in (\ref{1.8}), to
estimate $\|\Theta-\overline{\Theta}\|_{L^p(\R_+)}$, it suffices to
prove

$$ \|\Theta-\theta_+\|_{L^p(\R_+)}\to 0,\ \ \mathrm{as}\ \ \kappa\to 0,\ p\geq1.$$
Because if $M_0>0$ is a constant independent of $\alpha$,
$$\|\Theta-\theta_+\|^p_{L^p(\R_+)}\leq
M_0\|\Theta-\theta_+\|_{L^1(\R_+)},$$ the only thing we need to
proof is
$$\lim_{\kappa\to 0}\|\Theta-\theta_+\|_{L^1(\R_+)}=0.$$

 In fact we set $sgn_{\eta}(s)=\left\{
                   \begin{array}{ll}
                     1, & \hbox{$s> \eta$;} \\
                     s/\eta, & \hbox{$-\eta\leq s\leq\eta$;} \\
                     -1, & \hbox{$s< -\eta$.}
                   \end{array}
                 \right.
$, $I_{\eta}(s)=\int_0^s sgn_{\eta}(s) ds$ and $\eta>0$. Both side
of (\ref{1.11})$_1$ multiply by $sgn_{\eta}(\Theta-\theta_+)$ and
integrate in $(0,+\infty)\times (0,t)$ we can get
\begin{eqnarray*}
\int_0^t\left(\int_{0}^{+\infty}I_{\eta}(\Theta-\theta_+)dx\right)_{\tau}
d\tau&=&-a\int_0^t(\ln\Theta)_x(0,t)sgn_{\eta}(\Theta-\theta_+)(0,t)d\tau\\
&&\quad-a\int_0^t\int_{0}^{+\infty}(\ln\Theta)_x^2sgn'_{\eta}(\Theta-\theta_+)dxd\tau.
\end{eqnarray*}

When $\eta\to 0$ and use (\ref{2.7}) , (\ref{2.11})  we can get
\begin{eqnarray}\label{2.36}
&&\|\Theta-\theta_+\|_{L^1(\R_+)}
+a\int_0^t\int_0^{+\infty}(\ln\Theta)_x^2sgn'_{\eta}(\Theta-\theta_+)dxd\tau\nonumber\\
&&=a\int_0^t(\ln\Theta)_x(0,t)sgn_{\eta}(\Theta-\theta_-)(0,t)d\tau+\|\Theta_0-\theta_+\|_{L^1(\R_+)}.
\end{eqnarray}
Similar as (\ref{2.27}), when we integrate (\ref{2.25}) in
$\R_+\times (0,t)$ and combine with (\ref{2.7})
 we can get that there
exist constant $M_0>0$ independent of $\alpha$ such that
 \begin{equation}\label{2.37}
\|(\ln\Theta)_{xx}\|^2+a\int_0^t\|\partial_x^3(\ln\Theta)\|^2d\tau\leq
M_0\alpha^3+M_0a\alpha^3+M_0a^{-1}\alpha.
 \end{equation}
Since $a=\kappa p_+(\gamma-1)/(\gamma R^2)$  and $\kappa \to 0$, we
can choose $\alpha^{-1/2}=\kappa<1$, use (\ref{2.7}) and
(\ref{2.37}) such that
 (\ref{2.36}) is meant  if constant $M_0>0$ is independent of $\alpha$, we
 have
$$\|\Theta-\theta_+\|_{L^1(\R_+)}\leq M_0t((a\alpha)^{3/4}+a^{5/4}\alpha)+M_0\alpha^{-1}\leq M_0\kappa^{3/8}(t+1),$$
so we get $\|(V-\overline{V},\Theta-\overline{\Theta})\|_{L^p}\to 0$
as $\kappa\to 0$ with any $t\in [0,T]$.

 It remains to estimate
$\|U-\overline{U}\|_{L^p}$. To do so, both side of (\ref{2.24})
multiply by $sgn_{\eta}((\ln\Theta)_x)$ then integrate in
$\R_+\times(0,t)$ we can get
\begin{eqnarray*}&&\int_0^t\left(\int_{\R_+}I_{\eta}((\ln\Theta)_x)dx\right)_{\tau}d\tau+a\int_0^t\int_{\R_+}\frac{(\ln\Theta)^2_{xx}}{\Theta}sgn'_{\eta}((\ln\Theta)_x)dxd\tau\\
&&=-a\int_0^t\theta_-^{-1}(\ln\Theta)_{xx}(0,\tau)sgn_{\eta}((\ln\Theta)_x)(0,\tau)d\tau\\
&&=-\int_0^t\theta_-^{-1}\Theta_{\tau}(0,t)sgn_{\eta}((\ln\Theta)_x)(0,\tau)d\tau=0.\end{eqnarray*}
Again let $\eta\to 0$ we can get that there exist constant $M_0>0$
independent of $\alpha$ such that
\begin{eqnarray}&&\int_{\R_+}|(\ln\Theta)_x|dx+a\int_0^t\int_{\R_+}\frac{(\ln\Theta)^2_{xx}}{\Theta}sgn'_{\eta}((\ln\Theta)_x)dxd\tau\nonumber\\
&&\leq M_0.\label{2.38}\end{eqnarray}

Use the definition of $U$ in (\ref{1.10}) and combine with
(\ref{2.7}),(\ref{2.37}) and (\ref{2.38}) we know that
\begin{eqnarray*}
&&\|U-\overline{U}\|^p_{L^p}\leq
M_0\kappa^p\|(\ln\Theta)_x\|_{L^1}\|(\ln\Theta)_x\|^{(p-1)/2}\|(\ln\Theta)_{xx}\|^{(p-1)/2}\\
&&\leq M_0\kappa^p\alpha^{p-1}.
\end{eqnarray*}
Remind that $\alpha=\kappa^{-1/2}$, so we can get
$$\lim_{\kappa\to 0}\|U-\overline{U}\|_{L^p}=0.$$
The proof of Lemma \ref{lem2.4} is therefore complete, which also
means $(V,U,\Theta)$ is viscous contact discontinuity.$\Box$

\section{Reformulation and Main result}
Let $(v,u,\theta)$ be the solution of the problem (\ref{1.6}), and
let $(V,U,\Theta)$ be the viscous contact discontinuity constructed
in (\ref{1.12}). Denote
\begin{eqnarray}
&&\varphi(x,t)=v(x,t)-V(x,t),\nonumber\\
&&\psi(x,t)=u(x,t)-U(x,t),\nonumber\\
&&\zeta(x,t)=\theta(x,t)-\Theta(x,t).\label{1.14}
\end{eqnarray}
Combining (\ref{1.12}) and (\ref{1.6}), the original problem can be
reformulated as
\begin{equation}\left\{
    \begin{array}{lll}
      \varphi_t=\psi_x, &  \\
      \psi_t-(\frac{R\Theta}{vV}\varphi)_x+(\frac{R\zeta}{v})_x=-\mu(\frac{U_x}{vV}\varphi)_x+\mu(\frac{\psi_x}{v})_x-F, & \\
      \frac{R}{\gamma-1}\zeta_t+\frac{R\theta}{v}(\psi_x+U_x)-\frac{R\Theta}{V}U_x
=\kappa(\frac{\zeta_x}{v})_x-\kappa(\frac{\Theta_x\varphi}{vV})_x+\mu(\frac{{u_x}^2}{v}-\frac{{U_x}^2}{V})-G,
&\\
\zeta(0,t)=0,&\\
\left(\frac{R\theta_-}{V+\varphi}-\mu\frac{U_x+\psi_x}{V+\varphi}\right)(0,t)=p_+,&\\
 (\varphi,\psi,\zeta)(x,0)=(\varphi_0,\psi_0,\zeta_0)=(v_0-V_0,u_0-U_0,\theta_0-\Theta_0).&\\
\end{array}
  \right.
\label{1.15}\end{equation} From (\ref{1.12}) it is easy to check
that the initial-boundary data in (\ref{1.15}) satisfies the
compatible condition, and
$$(\varphi,\psi,\zeta)(x,0)=(\varphi_0,\psi_0,\zeta_0)\to (0,0,0)\ \ \mathrm{as}\ \ x\to +\infty.$$
To state our main result, we assume throughout of this section that
$$(\varphi_0,\zeta_0)(x)\in H_0^1(0,\infty),\ \ \psi_0(x)\in H^1(0,\infty).$$
Moreover, for an interval $I\in [0,\infty)$ , we define the function
space
$$X(I)=\left\{(\varphi,\psi,\zeta)\in C(I,H^1)|\varphi_x\in
L^2(I;L^2), (\psi_x,\zeta_x)\in L^2(I;H^1)\right\}.$$

Our  main results of this paper now reads as follows.
\begin{thm}\label{thm1.1}
There exist  positive constants $C>1$ and $\eta_0$ such that if
   $\|(v_0-V_0,
u_0-U_0, \theta_0-\Theta_0)\|_{L^2}\leq \eta_0,$ $\|(v_{0x}-V_{0x},
u_{0x}-U_{0x}, \theta_{0x}-\Theta_{0x})\|_{L^2}\leq C,$(\ref{1.15})
has a unique global solution $(\varphi,\psi,\zeta)$ satisfying
$(\varphi,\psi,\zeta)\in X([0,\infty))$ and
 $$\sup_{x\in \R_+}|(\varphi,\psi,\zeta)|\to 0,\ as\ t\to \infty.$$
\end{thm}

In this section, to study the asymptotic behavior of the solution to
the free boundary problem (\ref{1.6}), we will do some preparation
lemmas and  list some priori estimates   which are important to the
proof of Theorem \ref{thm1.1}.

We shall prove Theorem \ref{thm1.1} by combining the local existence
and the global-in-time  priori estimates. Since the local existence
of the solution  is well known (see, for example,\cite{HMS}), we
omit it here for brevity. to prove the global existence part of
Theorem \ref{thm1.1}, it is sufficient to establish the following
priori estimates.

\begin{pro}\label{pro2.2}{\rm(A priori estimate)} Let $(\varphi,\psi,\zeta)\in
X([0,t])$ be a solution of problem (\ref{1.15}) for some $t>0$. Then
there exist positive  constants $C(\delta_0)<1$ and $C$ which are
all independent of $t$ and $(v,\theta)$ , such that if $m\leq
v,\theta\leq M$
 and $1<N(t)=\sup_{0\leq \tau\leq t}\|(\varphi,\psi,\zeta)\|_1\leq
C$,
 it holds that
\begin{eqnarray}\label{2.1}
&&\sup_{0\leq \tau\leq t}\|(\psi,\varphi,\zeta)\|^2(t)+\int_0^t\|(\psi_x,\zeta_x)\|^2(\tau)d\tau\nonumber\\
&&\leq C\|(\varphi_0,\psi_0,\zeta_0)\|^{1/2}+C(\delta_0)
.\nonumber\\
&&\sup_{0\leq \tau\leq t}\|(\psi_x,\varphi_x,\zeta_x)\|^2(t)+\int_0^t\left(\|\varphi_x\|^2(\tau)+\|(\psi_x,\zeta_x)\|^2_1(\tau)\right)d\tau\nonumber\\
&&\leq C\|(\varphi_0,\psi_0,\zeta_0)\|_1+C(\delta_0).\end{eqnarray}
\end{pro}

\section{Proof of Theorem \ref{thm1.1}}

 Under the preparations in last section, the main task here is to finish (\ref{2.1}).  This part we also do some preparations. we must use the results
\begin{eqnarray}\label{3.1}
&&|V_x|\leq C|\Theta_x|,\nonumber\\
&&|\Theta_x|^2\leq
C\|(\ln\Theta)_x\|\|(\ln\Theta)_{xx}\|,\nonumber\\
&&|U_x|\leq C |(\ln\Theta)_{xx}|,\nonumber\\
 &&|U_x|^2\leq
C\|(\ln\Theta)_{xx}\|\|(\ln\Theta)_{xxx}\|,
\end{eqnarray}
 which follow from
(\ref{1.10})--(\ref{1.12}) . Also we set $C(\delta_0)$ stands for
small constants about  $\delta_0$, $\|(\varphi_0,\psi_0,\zeta_0)\|$
is asked suitably small, $C_v=\frac{R}{\gamma-1}$ and
$$\epsilon_1\ll\epsilon_2\ll \epsilon_3.$$

Before establishing (\ref{2.1}), we first estimate the value of
$\varphi(0,t)$ on the boundary $x=0$ by the boundary condition
(\ref{1.15}). Let $\varphi(t)=\varphi(0,t)$. Since
$U_x(0,t)=V_t(0,t)=0$, the boundary condition of (\ref{1.15}) yields
\begin{equation}\label{3.2}
\frac{R\theta_-}{v_-+\varphi(t)}-\mu\frac{\varphi_t(t)}{v_-+\varphi(t)}=p_+,\
t>0
\end{equation}
Direct computation gives
\begin{equation}\label{3.3}
\varphi_t(t)=-\frac{p_+}{\mu}\varphi(t),\ \varphi(0)=\varphi_0(0).
\end{equation}
It follows then that
\begin{equation}\label{3.4}
\varphi(t)=\varphi_0(0)e^{-p_+t/\mu}.
\end{equation}

\begin{lem}\label{lem3.2} If $C(\delta_0)>0$ is a small constant about
$\delta_0$
\begin{eqnarray*}&&\int_0^t\int_{\R_+}\Theta_x^2(\varphi^2+\zeta^2)dxd\tau\leq
C(\delta_0)\int_0^t\|(\varphi_x,\zeta_x)\|^2d\tau+C\varphi_0(0).
\end{eqnarray*}
\end{lem}
\pf Because if $x>0$
\begin{eqnarray*}&&\frac{\varphi^2}{x+1}=\int_0^x\big(\frac{2\varphi\varphi_x}{x+1}-\frac{\varphi^2}{(x+1)^2}\big)dx+\varphi^2(0,t)\\
&&=\int_0^x\left(\varphi_x^2-(\varphi_x-\frac{\varphi}{x+1})^2\right)dx+\varphi^2(0,t)\leq\int_0^x\varphi_x^2dx+C\varphi_0(0)e^{-p_+t/\mu}\leq
\|\varphi_x\|^2+C\varphi_0(0)e^{-p_+t/\mu}.\end{eqnarray*} similar
as above we can
get\begin{eqnarray*}&&\frac{\zeta^2}{x+1}\leq\int_0^x\zeta_x^2dx\leq
\|\zeta_x\|^2.\end{eqnarray*}

As to
\begin{eqnarray*}&&\int_0^t\int_0^{+\infty}\Theta_x^2(\varphi^2+\zeta^2)dxd\tau\\
&&\leq
\int_0^t\int_0^{+\infty}\Theta_x^2(x+1)\frac{(\varphi^2+\zeta^2)}{1+x}dxd\tau\\
&&\leq
\int_0^t\big(\int_0^{+\infty}\Theta_x^2(1+x)dx\big)\|(\varphi_x,\zeta_x)\|^2d\tau+C\varphi_0(0),
\end{eqnarray*}
 use (\ref{2.7})and (\ref{2.13}) we can get \begin{eqnarray*}&&\int_0^t\int_{\R_+}\Theta_x^2(\varphi^2+\zeta^2)dxd\tau\leq C(\delta_0)\int_0^t\|(\varphi_x,\zeta_x)\|^2d\tau+C\varphi_0(0),
\end{eqnarray*} and we finish this lemma.$\Box$

 Now, let's finish (\ref{2.1}) by the following lemmas.

\begin{lem}\label{lem3.1}
If  $\epsilon_1>0$  and  $C(\delta_0)>0$ are small constant about
$\delta_0$,  we can get
\begin{eqnarray*}
&&\int_{\R_+}(\varphi^2+\psi^2+\zeta^2)dx+\int_0^t\
\left\|(\psi_x,\zeta_x)\right\|^2\
d\tau\\
&&\leq C(\delta_0)+C\left\{\epsilon_1\int_0^t\ \|\varphi_x\|^2
d\tau+\|(\varphi_0,\psi_0,\zeta_0)\|^2\right\}.
\end{eqnarray*}
\end{lem}
\pf Set
$$\Phi(z)=z-\ln z-1,$$
$$\Psi(z)=z^{-1}+\ln z-1,$$ where $\Phi'(1)=\Phi(1)=0$ is a strictly
convex function around $z=1$. Similar to the proof in \cite{HMS},
 we deduce from (\ref{1.15}) that
\begin{eqnarray}\label{3.5}
&&\left(\frac{\psi^2}{2}+R\Theta\Phi\left(\frac{v}{V}\right)+C_v\Theta\Phi\left(\frac{\theta}{\Theta}\right)\right)_t\nonumber\\
&&\quad+\mu\frac{\Theta\psi_x^2}{v\theta}+\kappa\frac{\Theta\zeta_x^2}{v\theta^2}+H_x+Q
=\mu\left(\frac{\psi\psi_x}{v}\right)_x-F\psi-\frac{\zeta
G}{\theta},
\end{eqnarray}
where
$$H=R\frac{\zeta\psi}{v}-R\frac{\Theta\varphi\psi}{vV}+\mu\frac{U_x\varphi\psi}{vV}-\kappa\frac{\zeta\zeta_x}{v\theta}+\kappa\frac{\Theta_x\varphi\zeta}{v\theta V},$$
and
\begin{eqnarray*}
Q&=&p_+\Phi\left(\frac{V}{v}\right)U_x+\frac{p_+}{\gamma-1}\Phi\left(\frac{\Theta}{\theta}\right)U_x-\frac{\zeta}{\theta}(p_+-p)U_x-\mu\frac{U_x\varphi\psi_x}{vV}\\
&&\quad-\kappa
\frac{\Theta_x}{v\theta^2}\zeta\zeta_x-\kappa\frac{\Theta\Theta_x}{v\theta^2V}\varphi\zeta_x-2\mu\frac{U_x}{v\theta}\zeta\psi_x
+\kappa\frac{\Theta_x^2}{v\theta^2V}\varphi\zeta+\mu\frac{U_x^2}{v\theta
V}\varphi\zeta\\
&=:&\sum_{i=1}^9Q_i.
\end{eqnarray*}
Note that $p=R\theta/v$, $p_+=R\Theta/V$ and (\ref{1.10}), use
integrate by part and Cauchy-Schwarz inequality can get

\begin{eqnarray}\label{3.6}
Q_1+Q_2&=&Ra\left(\Phi\left(\frac{V}{v}\right)(\ln\Theta)_x\right)_x+\frac{Ra}{\gamma-1}\left(\Phi\left(\frac{\Theta}{\theta}\right)(\ln\Theta)_x\right)_x\nonumber\\
&&\quad-aR(\ln\Theta)_x\left(\frac{V\varphi_x\varphi-V_x\varphi^2}{Vv^2}\right)\nonumber\\
&&\quad-a\frac{p_+}{\gamma-1}(\ln\Theta)_x\left(\frac{\Theta\zeta_x\zeta-\Theta_x\zeta^2}{\Theta\theta^2}\right)\nonumber\\
&&\geq\left(p_+\Phi\left(\frac{V}{v}\right)U+\frac{p_+}{\gamma-1}\Phi\left(\frac{\Theta}{\theta}\right)U\right)_x\nonumber\\
&&\quad-\epsilon(\zeta_x^2+\varphi_x^2)-C\epsilon^{-1}\Theta_x^2(\zeta^2+\varphi^2).
\end{eqnarray}
Similarly, using $p-p_+=\frac{R\zeta-p_+\varphi}{v}$, we can get
\begin{equation}\label{3.7}
Q_3\geq\frac{R\zeta-p_+\varphi}{v}(\frac{\zeta}{\theta}U_x)\geq\left(\frac{R\zeta^2U}{v\theta}-\frac{p_+\zeta\varphi
U}{\theta
v}\right)_x-\epsilon(\zeta_x^2+\varphi_x^2)-C\epsilon^{-1}\Theta_x^2(\zeta^2+\varphi^2).
\end{equation}
And
\begin{eqnarray}\label{3.8}
(Q_4+Q_7)+(Q_5+Q_6+Q_8)+Q_9&\geq&
-C\epsilon^{-1}(\ln\Theta)_{xx}^2-\epsilon\psi_x^2\nonumber\\
&&-\epsilon\zeta_x^2-C\epsilon^{-1}\Theta_x^2(\zeta^2+\varphi^2)\nonumber\\
&&-C\epsilon^{-1}|(\ln\Theta)_{xx}|^2(\zeta^2+\varphi^2).
\end{eqnarray}
At the end we use the definition of $F$ and $G$ in (\ref{1.13}) then
combine with the general inequality skills as above to get
\begin{eqnarray}\label{3.9}
-F\psi-G\frac{\zeta}{\theta}&=&-\frac{\kappa a(\gamma-1)-\mu
p_+\gamma}{R\gamma}\left(\frac{(\ln\Theta)_{xx}}{\Theta}\right)_x\psi\nonumber\\
&&\quad+\frac{\mu p_+}{R\Theta}\left(\frac{\kappa
(\gamma-1)}{R\gamma}(\ln\Theta)_{xx}\right)^2\frac{\zeta}{\theta}\nonumber\\
&\leq&-\frac{\kappa a(\gamma-1)-\mu
p_+\gamma}{R\gamma}\left(\frac{(\ln\Theta)_{xx}}{\Theta}\psi\right)_x+\frac{\kappa
a(\gamma-1)-\mu
p_+\gamma}{R\gamma}\frac{(\ln\Theta)_{xx}}{\Theta}\psi_x\nonumber\\
&&\quad+\frac{\mu p_+}{R\Theta}\left(\frac{\kappa
(\gamma-1)}{R\gamma}(\ln\Theta)_{xx}\right)^2\frac{\zeta}{\theta}\nonumber\\
&\leq&-\frac{\kappa a(\gamma-1)-\mu
p_+\gamma}{R\gamma}\left(\frac{(\ln\Theta)_{xx}}{\Theta}\psi\right)_x+\epsilon\psi_x^2+C\epsilon^{-1}(\ln\Theta)_{xx}^2.
\end{eqnarray}
Integrating (\ref{3.6})$-$(\ref{3.9}) in $\R\times(0,t)$ , using
(\ref{2.7}), (\ref{2.11}) and the boundary condition about
$(\varphi,\psi,\zeta)$ of (\ref{1.15}), (\ref{3.4}) and
$\Theta_t(0,t)=a(\ln\Theta)_{xx}(0,t)=0$,
$\Phi\left(\frac{\Theta}{\theta}\right)(0,t)=0$ to estimate the
terms
$\left(p_+\Phi\left(\frac{V}{v}\right)U+\frac{p_+}{\gamma-1}\Phi\left(\frac{\Theta}{\theta}\right)U\right)_x$,
$\mu\left(\frac{\psi\psi_x}{v}\right)_x$,
$\left(\frac{(\ln\Theta)_{xx}\psi}{\Theta}\right)_x$ and $H_x$, we
can know that because
\begin{eqnarray*}
&&\int_0^t\left(Ra\Phi\left(\frac{V}{v}\right)|(\ln\Theta)_x|+\mu\left|\frac{\psi\psi_x}{v}\right|+\left|\frac{(\ln\Theta)_{xx}\psi}{\Theta}\right|+|H|\right)(0,\tau)d\tau\\
&&\leq
C\int_0^t|(\ln\Theta)_x|(0,\tau)\left(|\varphi|(\tau)+|\varphi_\tau|(\tau)\right)d\tau\leq
C(\delta_0),
\end{eqnarray*}
 it is easy to get
\begin{eqnarray}\label{3.10}
&&\left|\int_0^t\int_{\R_+}\left(\mu\left(p_+\Phi\left(\frac{V}{v}\right)U+\frac{p_+}{\gamma-1}\Phi\left(\frac{\Theta}{\theta}\right)U\right)_x+\left(\frac{(\ln\Theta)_{xx}\psi}{\Theta}\right)_x+H_x\right)dxd\tau\right|\nonumber\\
&&\leq C(\delta_0).
\end{eqnarray}

In the end from  combine with the estimates from $Q_1$ to $Q_9$,
(\ref{3.10}) and Lemma \ref{lem3.2} we have
\begin{eqnarray*}
&&\int_{\R_+}\left(R\theta\Phi\left(\frac{v}{V}\right)+\frac{1}{2}\psi^2+C_v\theta\Phi\left(\frac{\theta}{\Theta}\right)\right)dx+\int_0^t\
\left\|\left(\psi_x/(\sqrt{v\theta}),\zeta_x/(\theta\sqrt{v})\right)\right\|^2\
d\tau\\
&&\leq
C\epsilon_1^{-1}\int_0^t\int_0^{\infty}\Theta_x^2(\varphi^2+\zeta^2)\
dxd\tau+C\left\{\epsilon_1\int_0^t\ \|\varphi_x\|^2\
d\tau+\|(\varphi_0,\psi_0,\zeta_0)\|^2\right\}+C(\delta_0)\\
&&\leq
C(\delta_0)\int_0^t\|(\psi_x,\zeta_x)\|^2d\tau++C\left\{\epsilon_1\int_0^t\
\|\varphi_x\|^2\
d\tau+\|(\varphi_0,\psi_0,\zeta_0)\|^2\right\}+C(\delta_0).
\end{eqnarray*}
Use the condition of Proposition \ref{pro2.2} which is $m\leq
v,\theta\leq M$, we finish this lemma.$\Box$

\begin{lem}\label{lem3.3} For a small $\epsilon_2>0$, $C(\delta_0)>0$ is a small constant about
$\delta_0$,we have
\begin{eqnarray*}
&&\|(\varphi,\psi,\zeta)\|^2+\|(\psi_x,\zeta_x)\|^2+\int_0^t\|(\psi_{xx},\zeta_{xx})\|^2dt\\
&&\leq
C\left(\|(\psi_{0x},\zeta_{0x})\|^2+\epsilon_2^{-1}\|(\varphi_0,\psi_0,\zeta_0)\|^2\right)+C\epsilon_2^{-1}\int_0^t\|\varphi_x\|^2d\tau+C(\delta_0)+C\varphi_0(0).
\end{eqnarray*}.
\end{lem}

\pf First to get the estimate of $\|\psi_x(t)\|$ ,multiply both side
of  (\ref{1.15})$_2$ to $\psi_{xx}$ to get
\begin{eqnarray*}
&&\left(\frac{\psi_x^2}{2}\right)_t+\mu\frac{\psi_{xx}^2}{v}
=\mu\frac{\psi_x v_x}{v^2}\psi_{xx}+\mu\left(\frac{U_x\varphi}{v
V}\right)_x\psi_{xx}\\
&&\quad-R\left(\frac{\Theta\varphi}{v
V}\right)_x\psi_{xx}+R\left(\frac{\zeta}{v}\right)_x\psi_{xx}+F\psi_{xx}+(\psi_t\psi_x)_x:=\sum_{i=1}^6I_i.
\end{eqnarray*}
use last inequality integrate in $\R_+\times (0,t)$
\begin{eqnarray}
&&\|\psi_x(t)\|^2+\int_0^t\|\psi_{xx}(\tau)\|^2d\tau\nonumber\\
&&\quad\leq
C\|\psi_{0x}\|^2+C\sum_{i=1}^6\left|\int_0^t\int_0^\infty
I_idxd\tau\right|.\label{3.11}
\end{eqnarray}

Now deal with   $\iint|I_i|dxd\tau$ in the right side of
(\ref{3.11}). Using $\epsilon$ small and $v=\varphi+V$ ,
$R\Theta/V=p_+$ and (\ref{2.9}), (\ref{2.10}), (\ref{3.1})to get
\begin{eqnarray}
&&\int_0^t \int_0^{\infty}|I_1| dx d\tau\leq C\int_0^t
\int_0^{\infty}|V_x|
|\psi_{x}||\psi_{xx}|dxd\tau+C\int_0^t\int_0^{\infty}
|\varphi_x||\psi_{x}||\psi_{xx}|dx d\tau\nonumber\\
&&\quad\leq C\int_0^t\|V_x\|\|\psi_x\|_{L^\infty}\|\psi_{xx}\|
d\tau+C\int_0^t \|\psi_x\|_{L^\infty}\|\varphi_x\|\|\psi_{xx}\|
d\tau\nonumber\\
&&\quad\leq \epsilon\int_0^t\|\psi_{xx}\|^2d\tau+
C\epsilon^{-1}\int_0^t\|\psi_x\|^2 \|V_x\|^4d\tau+C\int_0^t
\|\psi_x\|^{1/2}\|\varphi_x\|\|\psi_{xx}\|^{3/2}
d\tau\nonumber\\
&&\quad\leq C\epsilon\int_0^t\|\psi_{xx}\|^2d\tau+
C(\delta_0)\int_0^t\|\psi_x\|^2
d\tau\nonumber\\
&&\quad\quad+C\epsilon^{-1}\sup_{t}\|\varphi_x\|^4\int_0^t
\|\psi_x\|^2 d\tau .\label{3.12}
\end{eqnarray}
Next we use the definition of $(V,U,\Theta)$ in
(\ref{1.10}),(\ref{1.12}), Cauchy-Schwarz inequality and
(\ref{2.11}), (\ref{2.12}), (\ref{3.1}) to get
\begin{eqnarray}
&&\int_0^t\int_0^\infty|I_2|dxd\tau\nonumber\\
&&\quad\leq
C\int_0^t\int_0^\infty\left(|U_{xx}||\varphi|+|U_x||\varphi_x|+|U_x||V_x||\varphi|+|U_x||\varphi||\varphi_x|\right)|\psi_{xx}|dxd\tau\nonumber\\
&&\quad\leq\epsilon\int_0^t\|\psi_{xx}\|^2d\tau+\frac{C}{\epsilon}\int_0^t\|\varphi\|_{L^\infty}^2\|U_{xx}\|^2d\tau
+\frac{C}{\epsilon}\int_0^t\|U_x\|_{L^\infty}^2\|\varphi_x\|^2d\tau\nonumber\\
&&\qquad+\frac{C}{\epsilon}\int_0^t\|\varphi\|_{L^\infty}^2\|V_x\|^2\|U_x\|_{L^\infty}^2d\tau
+\frac{C}{\epsilon}\int_0^t\|\varphi\|_{L^\infty}^2\|U_x\|_{L^\infty}^2\|\varphi_x\|^2d\tau\nonumber\\
&&\quad\leq\epsilon\int_0^t\|\psi_{xx}\|^2d\tau+C(\delta_0)+C(\delta_0)\int_0^t\|\varphi_x\|^2d\tau.\label{3.13}
\end{eqnarray}
The same as (\ref{3.12}) and (\ref{3.13}), we  the definition of $F$
in(\ref{1.13}) and (\ref{2.9})--(\ref{2.11}),(\ref{3.1})we can get
the estimates about $I_3$ to $I_5$ as following.
\begin{eqnarray}
&&\int_0^t\int_0^\infty(|I_3|+|I_4|+|I_5|)dxd\tau\nonumber\\
&&\leq
C\int_0^t\int_0^\infty\left(|\Theta_{x}||\varphi|+|\Theta||\varphi_x|+|\Theta||V_x||\varphi|
+|\Theta||\varphi||\varphi_x|\right)|\psi_{xx}|dxd\tau\nonumber\\
&&\quad+C\int_0^t\int_0^\infty\left(|\zeta_x|+|\zeta||V_x|+|\zeta||\varphi_x|\right)|\psi_{xx}|dxd\tau\nonumber\\
&&\quad+\epsilon\int_0^t\|\psi_{xx}\|^2d\tau+\frac{C}{\epsilon}\int_0^t\|F\|^2d\tau\nonumber\\
&&\leq\epsilon\int_0^t\|\psi_{xx}\|^2d\tau+C\epsilon^{-1}\int_0^t\|\varphi_x\|^2d\tau+\frac{C}{\epsilon}\int_0^t\|\varphi_x\|^2d\tau+\frac{C}{\epsilon}\int_0^t\int_0^\infty
V_x^2\varphi^2dx
d\tau\nonumber\\
&&\quad+\epsilon\int_0^t\|\psi_{xx}\|^2d\tau+\frac{C}{\epsilon}\int_0^t\int_0^\infty\left(\zeta_x^2+V_x^2\zeta^2\right)dxd\tau+\frac{C}{\epsilon}\sup_t\|(\varphi,\zeta)\|\|(\varphi_x,\zeta_x)\|\int_0^t\|\varphi_x\|^2d\tau\nonumber\\
&&\quad+\epsilon\int_0^t\|\psi_{xx}\|^2d\tau+C(\delta_0)+C\varphi_0(0).\label{3.14}
\end{eqnarray}

Because Lemma \ref{lem3.2} and Lemma \ref{lem3.1}, we know
$\|(\varphi,\zeta)\|$ is suitably small when $C(\delta_0)$ and
$\|(\varphi_0,\psi_0,\zeta_0)\|$ small. So there exist a small
constant $\delta$ about $\|(\varphi_0,\psi_0,\zeta_0)\|$ and
$\delta_0$  such that
$$\frac{C}{\epsilon}\sup_t\|(\varphi,\zeta)\|\|(\varphi_x,\zeta_x)\|\int_0^t\|\varphi_x\|^2d\tau\leq C\delta\int_0^t\|\varphi_x\|^2d\tau,$$
here $\|\varphi_x\|^2+\int_0^t\|\varphi_x\|^2d\tau\leq C$ can be
established in Lemma \ref{lem3.4}. Therefore
\begin{eqnarray*}
&&\int_0^t\int_0^\infty(|I_3|+|I_4|+|I_5|)dxd\tau\\
&&\leq
\epsilon\int_0^t\|\psi_{xx}\|^2d\tau+\frac{C}{\epsilon}\int_0^t\int_0^\infty\
(\zeta_x^2+\varphi_x^2)dxd\tau
+C(\delta_0)+C\varphi_0(0).\end{eqnarray*}

At last we integrate by part to the term about $I_6$. Because
$\psi_x=\varphi_t$,
\begin{equation}\label{3.15}
\psi_t\psi_x(0,t)=\psi_t\varphi_t(0,t)=(\psi\varphi_t)_t(0,t)-\varphi_0(0)\frac{p_+^2}{\mu^2}\psi(0,t)
e^{-p_+t/\mu}
\end{equation}
and $\psi(0,t)=C(\ln\Theta)_x(0,t)$,

we have
\begin{eqnarray}
&&\left|\int_0^t\int_0^\infty
I_6dxd\tau\right|=\left|\int_0^t(\psi_t\varphi_t)(0,\tau)d\tau\right|\nonumber\\
&&\leq
C|\psi(0,t)\varphi_t(0,t)|+C|\psi(0,0)\varphi_t(0,0)|+C\int_0^t\varphi_0(0)|\psi|(0,t)e^{-p_+\tau/\mu}d\tau\nonumber\\
&&\leq C(\delta_0).\label{3.16}
\end{eqnarray}

In all we can get
\begin{eqnarray}\label{3.17}
&&\int_0^t\ \int_0^{\infty}\ \sum_{i=1}^6|I_i|\ dx\ d\tau\nonumber\\
&&\leq C\int_0^t\
\epsilon\|\psi_{xx}\|^2\ d\tau+CN^4(t)\epsilon^{-1}\int_0^t\|\psi_x\|^2d\tau\nonumber\\
&&\quad+C\epsilon^{-1}\int_0^t\ \|(\varphi_x,\psi_x,\zeta_x)\|^2\
d\tau+C(\delta_0)+C\varphi_0(0).
\end{eqnarray}

So (\ref{3.11}) can be change to
\begin{eqnarray}
\|\psi_x(t)\|^2+\int_0^t\|\psi_{xx}(\tau)\|^2d\tau &\leq&
C(\delta_0)+C\varphi_0(0)+C\epsilon^{-1}\int_0^t\|(\varphi_x,\psi_x,\zeta_x)\|^2d\tau\nonumber\\
&&\quad+CN^4(t)\epsilon^{-1}\int_0^t\|\psi_x\|^2d\tau+C\|\psi_{0x}\|^2.
\label{3.18}
\end{eqnarray}

The estimate about  $\|\zeta_x\|$ is similar to $\|\psi_x\|$, use
(\ref{1.15})$_3$ multiply $\zeta_{xx}$ then integrate in
$Q_t=\R_+\times(0,t)$ to get
\begin{eqnarray}
&&\|\zeta_x\|^2+\int_0^t
\|\zeta_{xx}\|^2\ d\tau\nonumber\\
&&\leq C \|\zeta_{0x}\|^2+
C\epsilon^{-1}\int_0^t\int_0^\infty\left(\psi_x^2+\zeta\psi_x^2+\zeta^2U_x^2+U_x^2\varphi^2\right)dxd\tau\nonumber\\
&&
\quad+C\int_0^t\int_{\R_+}|\zeta_x|(|\varphi_x|+|V_x|)|\zeta_{xx}|dxd\tau+C\epsilon^{-1}\int_0^t\int_0^\infty\left|\left(\frac{\Theta_x\varphi}{vV}\right)_x\right|^2dxd\tau
\nonumber\\
&&\quad+C\epsilon^{-1}\int_0^\infty\int_0^\infty(U_x^4+\psi_x^4)dxd\tau+C\epsilon^{-1}\int_0^t\|G\|^2d\tau\nonumber\\
&&=:C\|\zeta_{0x}\|^2+\sum_{i=1}^5J_i.\label{3.19}
\end{eqnarray}
Use the same method as (\ref{3.12})--(\ref{3.16})
$$
J_1\leq
C\epsilon^{-1}\int_0^t\|\psi_x\|^2d\tau+C\epsilon^{-1}N^2(t)\int_0^t\|U_x\|^2d\tau
\leq C\epsilon^{-1}\int_0^t\|\psi_x\|^2d\tau+C(\delta_0).
$$
Again use the same method as(\ref{3.12})--(\ref{3.16})
\begin{eqnarray*}
J_2&\leq& C\int_0^t\|\zeta_x\|_{L^\infty}\|\varphi_x\|\|\zeta_{xx}\|d\tau+C\int_0^t\|V_x\|\|\zeta_x\|_{L^\infty}\|\zeta_{xx}\|d\tau\\
&\leq&C\int_0^t\|\zeta_x\|^{1/2}\|\zeta_{xx}\|^{3/2}\|\varphi_x\|d\tau
+\epsilon\int_0^t\|\zeta_{xx}\|^2d\tau+C(\delta_0)\int_0^t\|\zeta_x\|^2d\tau\\
&\leq&2\epsilon\int_0^t\|\zeta_{xx}\|^2d\tau+C(\delta_0)\int_0^t\|\zeta_x\|^2d\tau
+C\epsilon^{-1}\sup_t\|\varphi_x\|^4\int_0^t\|\zeta_x\|^2d\tau.
\end{eqnarray*}

Because
\begin{eqnarray*}&&\left|\left(\frac{\Theta_x\varphi}{vV}\right)_x\right|^2\\
&&=|\frac{\Theta_{xx}\varphi}{vV}+\frac{\Theta_x\varphi_x}{vV}+\frac{\Theta_x\varphi}{vV}(-\frac{V_x+\varphi_x}{v^2}-\frac{V_x}{V^2})|^2\\
&&\leq
C\Theta_{xx}^2\varphi^2+C\Theta_x^2\varphi_x^2+C\Theta_x^2V_x^2\varphi^2+C\Theta_x^2\varphi^2\varphi_x^2,
\end{eqnarray*}
combine with $R\Theta/V=p_+$, use the same method as
(\ref{3.12})--(\ref{3.16}) to get
\begin{eqnarray*}
J_3&\leq&
C\epsilon^{-1}\int_0^t\|\varphi\|_{L^\infty}^2\|\Theta_{xx}\|^2d\tau
+C\epsilon^{-1}\int_0^t\|\Theta_x\|_{L^\infty}^2\|\varphi_x\|^2d\tau\\
&&\quad+C\epsilon^{-1}\int_0^t\ \int_0^{\infty}\ \Theta_x^2V_x^2\varphi^2\ dx\ d\tau\nonumber\\
&\leq&C(\delta_0)\int_0^t\|\varphi_{x}\|^2d\tau+C(\delta_0)+C\varphi_0(0)
.
\end{eqnarray*}
Use the definition $U$ and similar as (\ref{3.12}) (\ref{3.13}) that
we combine with Lemma \ref{lem2.3} to get
\begin{eqnarray*} J_4&\leq&
C(\delta_0)+C\epsilon^{-1}\int_0^t\|\psi_{x}\|_{L^\infty}^2\|\psi_x\|^2d\tau\\
&\leq& C(\delta_0)+C\epsilon^{-1}\int_0^t\|\psi_x\|^3\|\psi_{xx}\|d\tau\\
&\leq&
C(\delta_0)+C\int_0^t\left(\epsilon^{-2}\|\psi_x\|^2\|\psi_x\|^4+\epsilon^{2}\|\psi_{xx}\|^2\right)d\tau.
\end{eqnarray*}
Use the definition $G$ in (\ref{1.13}) combine with Lemma
\ref{lem2.3}
\begin{eqnarray*}
J_5=C\epsilon^{-1}\int_0^t\|G\|^2d\tau\leq C(\delta_0).
\end{eqnarray*}
Use the results from $J_1$ to $J_5$, the inequality (\ref{3.19}) can
be change to
\begin{eqnarray}
&&\|\zeta_x\|^2+\int_0^t\ \|\zeta_{xx}\|^2\ d\tau\nonumber\\
&&\leq
C\|\zeta_{0x}\|^2+C(\epsilon^{-3}+N^4(t))\int_0^t\|(\psi_x,\zeta_x)\|^2d\tau+C(\delta_0)\int_0^t\|\varphi_x\|^2d\tau+C(\delta_0)\nonumber\\
&&\quad+C\epsilon\int_0^t\|\psi_{xx}\|^2d\tau+C\varphi_0(0).\label{3.20}
\end{eqnarray}

 In fact when combine  with  (\ref{3.18}) and (\ref{3.20}), it is easy to get
\begin{eqnarray*}
&&\|(\varphi,\psi,\zeta)\|^2+\|(\psi_x,\zeta_x)\|^2+\int_0^t\|(\psi_{xx},\zeta_{xx})\|^2d\tau\\
&&\leq
C\left(\|(\psi_{0x},\zeta_{0x})\|^2+\epsilon^{-3}\|(\varphi_0,\psi_0,\zeta_0)\|^2\right)+C\epsilon^{-1}\int_0^t\|\varphi_x\|^2d\tau+C(\delta_0)+C\varphi_0(0).
\end{eqnarray*}
$\Box$

\begin{lem}\label{lem3.4} For a small $\epsilon_3>0$  and $C(\delta_0)>0$ is a small constant about
$\delta_0$  , $C(\delta_0)>0$ is a small constant about $\delta_0$,
 we can get
\begin{eqnarray}
\|\varphi_x\|^2+\int_0^t\|\varphi_x\|^2 d\tau&\leq&
C\|\varphi_{0x}\|^2+C\epsilon_3^{-1}\|(\varphi_0,\psi_0,\zeta_0)\|^2
+C(\delta_0)+C\varphi_0(0) .\label{3.21}
\end{eqnarray}
\end{lem}
\pf Set $\bar{v}=\frac{v}{V}$ take it into (\ref{1.15})$_1$,
(\ref{1.15})$_2$  ($p=R \theta/v$) to get
$$
\psi_t+p_x=\mu\left (\frac{\bar{v}_x}{\bar{v}}\right)_t-F,
$$
Both sides of last equation multiply $\bar v_x/\bar v$ to get
\begin{eqnarray}
&&\left(\frac{\mu}{2}\left (\frac{\bar v_x}{\bar
v}\right)^2-\psi\frac{\bar v_x}{\bar v}\right)_t
+\frac{R\theta}{v}\left (\frac{\bar v_x}{\bar
v}\right)^2+\left(\psi\frac{\bar v_t}{\bar v}\right)_x\nonumber\\
&&\quad=\frac{\psi_x^2}{v}+U_x\left
(\frac{1}{v}-\frac{1}{V}\right)\psi_x+\frac{R\zeta_x}{v}\frac{\bar
v_x}{\bar v}-\frac{R\theta}{v}\left
(\frac{1}{\Theta}-\frac{1}{\theta}\right)\Theta_x\frac{\bar
v_x}{\bar v}+F\frac{\bar v_x}{\bar v}.\label{3.22}
\end{eqnarray}
On the boundary condition we have
$$\psi(0,\tau)\frac{{\bar v}_\tau(0,\tau)}{{\bar v}(0,\tau)}=\psi(0,\tau)\frac{\varphi_\tau(0,\tau)}{v(0,\tau)}$$
 When use the definition $\psi(0,t)$ in (\ref{1.15}) and $\varphi_t(0,t)$ in (\ref{3.3}) we can get
 \begin{equation}\label{3.23}
 \int_0^t\left|\psi(0,\tau)\frac{{\bar v}_\tau(0,\tau)}{{\bar
 v}(0,\tau)}\right|d\tau\leq C(\delta_0)
 \end{equation}
On the other hand if we integrate (\ref{3.22}) in $R_+\times(0,t)$,
 (\ref{3.22}) is changed to
\begin{eqnarray}\label{3.24}
&&\int_{\R_+}\left(\frac{\mu}{2}\left (\frac{\bar v_x}{\bar
v}\right)^2-\psi\frac{\bar v_x}{\bar
v}\right)dx-\int_{\R_+}\left(\frac{\mu}{2}\left (\frac{\bar
v_{x}(x,0)}{\bar v(x,0)}\right)^2-\psi_0\frac{\bar v_{x}(x,0)}{\bar v(x,0)}\right)dx\nonumber\\
&&\quad+\int_0^t\int_{\R_+}\left(\frac{R\theta}{v}\left (\frac{\bar
v_x}{\bar v}\right)^2+\left(\psi\frac{\bar v_t}{\bar
v}\right)_x\right)dxd\tau\nonumber\\
&&\leq C\epsilon^{-1}\left(\int_0^t\ \|(\zeta_x,\psi_x)\|^2\
d\tau+\int_0^t\ \int_0^{\infty}\ \Theta_x^2(\varphi^2+\zeta^2)\ dx\
d\tau\right)\nonumber\\
&&+C\epsilon^{-1}\int_0^t\int_0^{+\infty}U_x^2\varphi^2
dxd\tau+C\epsilon^{-1}\int_0^t\int_0^{+\infty}|F|^2dxd\tau+\epsilon\int_0^t\|\frac{\bar{v}_x}{\bar{v}}\|^2
d\tau+C(\delta_0).\end{eqnarray} Insert (\ref{3.23}) into
(\ref{3.24}) and combine with the definition of $F$ and $U$ and
Lemma \ref{lem2.3}, we can get
\begin{eqnarray}\label{3.25}&&\int_0^t\|\frac{\bar{v}_x}{\bar{v}}\|^2\
d\tau+\|\frac{\bar{v}_x}{\bar{v}}\|^2-C\epsilon^{-1}\|\psi\|^2-C\|\psi_0\|^2-C\int_0^{\infty}\frac{\bar{v}_x}{\bar{v}}(x,0)^2\
dx\nonumber\\ &&\leq C\epsilon^{-1}\left(\int_0^t\
\|(\zeta_x,\psi_x)\|^2\ d\tau+\int_0^t\ \int_0^{\infty}\
\Theta_x^2(\varphi^2+\zeta^2)\ dx\
d\tau+C(\delta_0)\right)\nonumber\\
&&\quad+\epsilon\int_0^t\|\frac{\bar{v}_x}{\bar{v}}\|^2
d\tau+C\|\varphi_{0x}\|^2.
\end{eqnarray}
Because $C_1(\varphi_x^2)-C_2V_x^2\leq
(\frac{\bar{v}_x}{\bar{v}})^2\leq C_3\varphi_x^2+C_4V_x^2$
($C_1,C_2,C_3,C_4$ stands for constants about $v$), combine with
Lemma \ref{lem3.2}-- \ref{lem3.1}  we change (\ref{3.25}) to
\begin{eqnarray}\label{3.26}\int_0^t\|\varphi_x\|^2\
d\tau+\|\varphi_x\|^2&\leq&
C\|\varphi_{0x}\|^2+\int_0^tC\epsilon_3^{-1}\|(\psi_x,\zeta_x)\|^2d\tau+C(\delta_0)+C\varphi_0(0).
\end{eqnarray}
 So we finish this lemma.$\Box$

 From Lemma \ref{lem3.1} to Lemma \ref{lem3.4} we know when $\delta_0$ and $\|(\varphi_0,\psi_0,\zeta_0)\|$ suitably small there exist a suitably small
 constant $\delta$
such that
$$\|(\varphi,\psi,\zeta)\|^2+\int_0^t\|(\psi_x,\zeta_x)\|^2d\tau\leq
 C\delta,$$ and $$\|(\varphi_x,\psi_x,\zeta_x)\|^2+\int_0^t\|(\psi_{xx},\zeta_{xx})\|^2\leq
 C.$$ Then we can get $C_5\leq |v|\leq
 C_6$ and $C_7\leq |\theta|\leq C_8$ when $\delta$ small, here $C_5$, $C_6,$ $C_7$ and
 $C_8$ are constants independent of $v$ and $\theta$ .
   So we can get
 (\ref{2.1}) in
 Proposition 2.2 .

To finish Theorem \ref{thm1.1} now we will proof $\sup_{x\in
\R_+}|(\varphi,\psi,\zeta)|\to 0,\ as\ t\to \infty.$

Because $\int_0^{\infty}\partial_x$(\ref{1.15})$_1\times 2\varphi_x\
dx$
 equals to
\begin{equation}
0=2\int_0^\infty\varphi_x\psi_{xx}dx-\frac{d}{dt}\|\varphi_x\|^2,\label{3.27}
\end{equation}
use  Cauchy-Schwarz inequality we get
$$2\int_0^\infty\varphi_x\psi_{xx}dx \leq C\left(\|\varphi_x\|^2+\|\psi_{xx}\|^2\right),$$

Again using Lemma \ref{lem3.3}--\ref{lem3.4} ,  we get
\begin{eqnarray}
&&\int_0^\infty\left|\frac{d}{dt}\|\varphi_x(t)\|^2\right|dt\nonumber\\
&&\leq
C\int_0^\infty\left(\|\varphi_x\|^2+\|\psi_{xx}\|^2\right)dt\nonumber\\
&&\leq
C\left(C(\delta_0)+\|(\varphi_0,\psi_0,\zeta_0)\|_1^2+C\varphi_0(0)\right).\label{3.28}
\end{eqnarray}
Similar as above, from Lemma \ref{lem3.1}$-$\ref{lem3.4} and combine
with Sobolev inequality we get
\begin{equation}
\int_0^\infty\left(\left|\frac{d}{dt}\|\psi_x(t)\|^2\right|+\left|\frac{d}{dt}\|\zeta_x(t)\|^2\right|\right)d\tau\leq
C\left(C(\delta_0)+\|(\varphi_0,\psi_0,\zeta_0)\|_1^2+C\varphi_0(0)\right).\label{3.29}
\end{equation}

It means
$$
\|(\varphi,\psi,\zeta)(t)\|_{L^\infty}^2\leq
2\|(\varphi,\psi,\zeta)(t)\|\|(\varphi_x,\psi_x,\zeta_x)(t)\|\to0\quad\mbox{when}\quad
t\to\infty.
$$
So we finish Theorem \ref{thm1.1}.$\Box$


\begin{thebibliography}{99}
\bibitem{X}
Z.P. Xin, On nonlinear stability of contact discontinuities.
Proceeding of 5th International Conferences on Hyperbolic Problems:
Theory, Numerics and Applications. Ed. Glimm, etc., World Sci.
Publishing, River Edge, NJ, 1996.
\bibitem{LX}
T.P. Liu, Z.P. Xin, Pointwise decay to contact discontinuities for
systems of viscous conservation laws. Asian J. Math., 1 (1997)
34-84.


\bibitem{HMS}
F.M. Huang, A. Matsumura, X. Shi, On the stability of contact
discontinuity for compressible Navier-Stokes equations with free
boundary. Osaka J. Math., 41 (1) (2004) 193-210.



\bibitem{HMX}
F.M. Huang, A. Matsumura, Z.P.Xin, Stability of contact
discontinuity for the 1-D compressible Navier-Stokes Equations.
Arch. Ratonal Mech. Anal., 179 (2005) 55-77.

\bibitem{HXY}
F.M. Huang, Z.P. Xin, T.Yang, Contact discontinuity with general
perturbations for gas motions. Adv. Math., 219 (2008) 1246-1297.

\bibitem{HZ}
F.M. Huang, H.J. Zhao, On the global stability of contact
discontinuity for compressible Navier-Stokes equations. Rend. Sem.
Mat. Univ. Padova, 109 (2003) 283-305.
\bibitem{HLM2009}
F.M. Huang, J. Li, A. Matsumura, Asymptotic stability of combination
of viscous contact wave with rarefaction waves for one-dimensional
compressible Navier-Stokes system. Arch.Ration.Mech.Anal.,
197(1)(2010)  89-116 .


\bibitem{MaSX} Ma S. X. ,Zero dissipation limit to strong  contact
discontinuity for the  $1-D$ compressible Navier-Stokes equations,
JDE., 248(2010) 95--110.
\bibitem{HH} Hakho Hong, Global stabillity of viscous contact wave
for $1-D$ compressible Navier-Stokes equations, JDE., 252(2012)
3482--3505.



\bibitem{Ma2001}
A. Matsumura, Inflow and outflow problems in the half space for a
one-dimensional isentropic model system of compressible viscous gas.
Meth. Appl. Anal., 8(4) (2001) 645-666.








\bibitem{TJJ2011}
T.Zheng, JW. Zhang, JN.Zhao, Asymptotic stability of viscous contact
discontinuity to an inflow problem for compressible Navier¨CStokes
equations,Nonlinear Analysis: Nonlinear Anal. Theor., 74(17)(2011)
6617-6639.

\bibitem{TZ2012}
T.Zheng and J.Zhao, On the stability of contact discontinuity for
Cauchy problem of compress Navier-Stokes equations with general
initial data. Sci. China. Math., 55(10)(2012) 2005-2026.



\bibitem{NYZ}
K. Nishihara, T. Yang, H.J. Zhao, Nonlinear stability of strong
rarefaction waves for compressible Navier-Stokes equations. SIAM J.
Math. Anal., 35 (2004) 1561-1597.

\bibitem{Z2013}
Tao Wang;Huijiang Zhao;Qingyang Zou,One-dimensional compressible
Navier-Stokes equations with large density oscillation,Kinetic and
Related Models, 6(3)(2013)649-670.

\bibitem{KM}
S. Kawashima, A. Matsumura, Asymptotic stability of traveling wave
solutions of systems for one-dimensional gas motion. Commun. Math.
Phys., 101 (1985) 97-127.
\bibitem{KMN}
S. Kawashima, A. Matsumura, K. Nishihara, Asymptotic behavior of
solutions for the equations of a viscous heat-conductive gas. Proc.
Japan Acad. Ser. A, 62 (1986) 249-252.
\bibitem{MN2001}
A. Matsumura, K. Nishihara, Large-time behaviors of solutions to an
inflow problem in the half space for a one-dimensional system of
compressible viscous gas, Comm.Math.Phys.222(2001)449-474.
\bibitem{KN}
S. Kawashima, S. Nishibata, P.C. Zhu, Asymptotic stability of the
stationary solution to the compressible Navier-Stokes equations in
the half space. Comm. Math. Phys., 240 (2003) 483-500.

\bibitem{KZ}
S. Kawashima, P.C. Zhu, Asymptotic stability of nonlinear wave for
the compressible Navier-Stokes equations in the half space. J. Diff.
Eqn., 244 (2008), 3151--3179.

\bibitem{L}
T.P. Liu, Nonlinear stability of shock waves for viscous
conservation laws, Mem. Amer.Math.Soc.,56(328),(1985).


\bibitem{MN1}
A. Matsumura, K. Nishihara, On the stability of travelling wave
solutions of a one-dimensional model system for compressible viscous
gas. Japan J. Appl. Math., 2 (1985) 17-25.

\bibitem{MN2}
A. Matsumura, K. Nishihara, Asymptotics toward the rarefaction waves
of the solutions of a one-dimensional model system for compressible
viscous gas. Japan J. Appl. Math., 3 (1985) 1-13.

\bibitem{MN3}
A. Matsumura, K. Nishihara, Global stability of the rarefaction wave
of a one-dimensional model system for compressible viscous gas.
Commun. Math. Phys., 165 (1992) 325-335.


\bibitem{SX}
Szepessy,A., Xin, Z.P., Nnlinear stability of viscous shck waves,
Arch.Ration.Mech.Anal.,122(1993)53-103.


\bibitem{SZ}
Szepessy,A., Zumbrun,K., Stability of rarefaction waves in viscous
media. Arch.Ration.Mech.Anal., 133(1996)249-298.




\end{thebibliography}
\end{document}